\documentclass[smallextended]{svjour3}       % onecolumn (second format)
\smartqed  % flush right qed marks, e.g. at end of proof
\usepackage{latexsym}
\usepackage{amsfonts,amsmath}
\newcommand{\re}{\text{\rm Re }}
\newcommand{\im}{\text{\rm Im }}

\begin{document}

\title{Approximation of sgn$(x)$ on two symmetric intervals by rational functions with fixed poles\thanks{Research is supported by RFBR-TUBITAK project No.14-01-91370/113F369}
}
\subtitle{Dedicated to Professor E.B. Saff on the occasion of his 70th birthday}

\titlerunning{Rational approximation of sgn $(x)$}        % if too long for running head

\author{Alexey  Lukashov         \and
       Dmitri Prokhorov %etc.
}


\institute{A. Lukashov \at
             Department of Mathematics \\
              Fatih University\\
              Turkey\\
              \email{alexey.lukashov@gmail.com}           %  \\
           \and
           A. Lukashov, D. Prokhorov \at
             Department of Mechanics and Mathematics\\ Saratov State University\\ Russian Federation\\
             \email{LukashovAL@info.sgu.ru, prokhorovdv@info.sgu.ru}}

\date{Received: date / Accepted: date}

\maketitle

\begin{abstract}
 Recently A. Eremenko and P. Yuditskii found explicit solutions of the best polynomial approximation problems of sgn(x) over two intervals in terms of conformal mappings onto special comb domains. We give analogous solutions for the best approximation problems of sgn(x) over two symmetric intervals by odd rational functions with fixed poles. Here the existence of the related conformal mapping is proved by using convexity of the comb domains along the imaginary axis.
 
\end{abstract}

\section{Introduction}
\label{intro}
Conformal mappings were used in solutions of extremal problems of approximation theory from the very beginning. In fact P.L. Chebyshev gave polynomials of the least deviation from zero on $[-1,1]$ in terms of the conformal mapping of $\mathbb{C}\backslash[-1,1]$ onto $\{|w|<1\}.$ Later N.I. Akhiezer \cite{akh} used conformal mappings of $\mathbb{C}\backslash([-1,a]\cup[b,1])$ onto a rectangle and of $\mathbb{C}\backslash([-1,a]\cup[c,d]\cup[b,1])$ onto a fundamental domain of a Schottky group to find the least deviating polynomials on two intervals $[-1,a]\cup[b,1].$ For several intervals generalizations of Akhiezer's constructions can be found in \cite{apt},\cite{bog},\cite{leb},\cite{luk1},\cite{luk2},\cite{lukpeh},\cite{luktys},\cite{peh},\cite{rob},\cite{sodyud}, \cite{wid} and references therein.

Recently one of those (comb domains) was heavily used to solve several extremal problems of the approximation theory by polynomials (see, for example, survey \cite{ereyudsurv}). In particular, the paper \cite{ereyud1} contains an explicit formula for the best approximation polynomials of the function sgn $(x)$ on two symmetric intervals in terms of a related comb function (the problem has many applications, see, for example \cite{SIAM1},\cite{SIAM2}).

The main goal of the paper is to show how more complicated comb domains can be used to solve extremal problems of the approximation theory by rational functions with fixed poles. Note also that additional difficulties appear here because of the question about the existence for desired conformal mappings.

Let $0<x_1<\dots<x_q<a<1<x_{q+1}<\dots<x_p$. For $B_0,B_1,\dots,B_p\in(0,\infty)$, ${\bf
B}=(B_0,B_1,\dots,B_p)$, $q\leq p$, $k_0,k_1,\dots,k_p,m\in\mathbb N$, ${\bf
k}=(k_0,k_1,\dots,k_p)$, denote $\Omega^{\bf k}_{q,m,p}({\bf B})$ a subdomain of the half-strip
$$\left\{w=u+iv:\;v>0,\;\;0<u<\pi\left(\sum_{j=0}^pk_j+m\right)\right\}$$ which is obtained after
deletion of the subset
$$E=\left\{w=u+iv:\;|u-\pi\sum_{j=0}^qk_j|\leq\arccos\left(\frac{\cosh B_0}{\cosh v}\right),
\;\;v\geq B_0\right\}$$ and the rays $$l_q=\{w=u+iv:\;u=\pi k_q,\;\;v\geq B_q\},$$
$$l_{q-1}=\{w=u+iv:\;u=\pi(k_q+k_{q-1}),\;\;v\geq B_{q-1}\},\dots,$$
$$l_1=\{w=u+iv:\;u=\pi(k_q+k_{q-1}+\dots+k_1),\;\;v\geq B_1\},$$ $$l_p=\left\{w=u+iv:\;
u=\pi\left(\sum_{j=0}^qk_j+m\right),\;\;v\geq B_p\right\},$$
$$l_{p-1}=\left\{w=u+iv:\;u=\pi\left(\sum_{j=0}^qk_j+m+k_p\right),\;\;v\geq
B_{p-1}\right\},\dots,$$
$$l_{q+1}=\left\{w=u+iv:\;u=\pi\left(\sum_{j=0,j\neq q+1}^pk_j+m\right),\;\;v\geq B_{q+1}\right\}.$$

Let $\phi(z)=\phi(z;q,m,p,{\bf k},{\bf B})$ be a conformal map from the first quadrant
$\{z=x+iy:\;\re z>0,\;\im z>0\}$ onto $\Omega^{\bf k}_{q,m,p}({\bf B})$ such that $\phi(a)=0$,
$\phi(1)=\pi(\sum_{j=0}^pk_j+m)$, $\phi(\infty)$ corresponds to $\infty$ reached along the
left-hand side of $l_p$.

The first theorem answers the question about the existence of desired conformal mapping.

\begin{theorem}
There is a unique vector ${\bf B^*}={\bf B^*}(x_1,\dots,x_p;a)$ such that $\phi(x_j)$ correspond to
$\infty$ reached along the left-hand side of $l_j$, $j=1,\dots,q$, $\phi(0)$ corresponds to
$\infty$ reached along the right-hand side of $l_1$, $\phi(x_j)$ correspond to $\infty$ reached
along the right-hand side of $l_j$, $j=q+1,\dots,p$.
\end{theorem}

 Theorem 1 conditions for conformal mapping $w=f(z)$ combine requirements of different kinds. In
the preimage domain of $z$-plane we have fixed points corresponding to the angular image points and
infinity. In the image domain of $w$-plane we have prescribed distances between vertical boundary
rays and the analytic description of the curvilinear boundary part up to a free parameter $B_0$ but
there is no information on ray tips.

Second theorem gives exact solution of an approximation problem in terms of the conformal mapping $\phi(z;q,m,p,{\bf k},{\bf B^*}).$

\begin{theorem}
The error of the best uniform approximation of the function $sgn (x),$ $|x|\in[a,1],$ by rational functions of the form
\begin{equation}
\label{form}
\frac{c_0x^n+c_{n-1}x^{n-1}+\ldots+c_0}{x^{2k_0-1}(x^2-x_1^2)^{k_1}\ldots(x^2-x_p^2)^{k_p}},\quad n=2(m-1+\sum_{j=0}^pk_j),
\end{equation}
is equal to
\begin{equation}	
\label{deviation}
L=\frac{1}{\cosh B_0^*(x_1,\ldots,x_p;a)},
\end{equation}
and the extremal rational function is given by $$f(x)=1-(-1)^{\sum_{j=0}^qk_j}L\cos\phi(x;q,m,p,\mathbf{k};\mathbf{B}^*),\quad x>0.$$
\end{theorem}

Other forms of approximating rational functions are considered analogously. For example, in the case of an even order pole at the origin instead of the left part of the boundary of the set $E$ the vertical half-line $$\left\lbrace w=u+iv: u=\pi\sum_{j=0}^qk_j,v\geq B_0\right\rbrace,$$ should be considered. Note also that for $p=0$ Theorem 2 was proved in \cite{ereyud2}.

\section{Proofs}
\label{sec:1}
\textit{Proof of Theorem 1}

First, map the upper half-disk $\{\zeta:\;|\zeta|<1,\;\im\zeta>0\}$ onto the first quadrant by a
fixed function $z=z(\zeta)$ so that $z(-1)=a$, $z(1)=1$ and $z(i)=0$. Denote
$z^{-1}(x_j):=e^{i\alpha_j}$, $j=1,\dots,p$; $z^{-1}(\infty):=e^{i\alpha_{p+1}}$,
$0<\alpha_{q+1}<\dots<\alpha_p<\alpha_{p+1}<\alpha_0<\alpha_1<\dots<\alpha_q<\pi$, where
$\alpha_0:={\pi\over2}$. Then $f(\zeta):=\phi(z(\zeta))$ maps the upper half-disk onto $\Omega^{\bf
k}_{q,m,p}({\bf B})$ so that $f(-1)=0$, $f(1)=\pi(\sum_{j=0}^pk_j+m)$, $f(e^{i\alpha_0})$
corresponds to $\infty$ reached along the right-hand side of $l_1$, $f(e^{i\alpha_{p+1}})$
corresponds to $\infty$ reached along the left-hand side of $l_p$, $f(e^{i\alpha_j})$ correspond to
$\infty$ reached along the left-hand side of $l_j$, $j=1,\dots,q$; $f(e^{i\alpha_j})$ correspond to
$\infty$ reached along the right-hand side of $l_j$, $j=q+1,\dots,p$.

The problem admits an equivalent posing for $f(\zeta)$ with the boundary conditions being
transformed from $\phi(z)$ to $f(\zeta)$.

Second, extend $f$ continuously onto the boundary of the upper half-disk, and apply the Schwarz
symmetry principle to extend $f$ onto the lower half-disk $\{\zeta:\;|\zeta|<1,\;\im\zeta<0\}$. The
extended function denoted also by $f$ maps the unit disk $\mathbb D$ conformally onto the domain
$\Omega:=\Omega^{\bf k}_{q,m,p}({\bf B})\cup(0,\pi(\sum_{j=0}^pk_j+m))\cup\overline{\Omega^{\bf
k}_{q,m,p}({\bf B})}$, where $\overline{\Omega^{\bf k}_{q,m,p}({\bf B})}$ is a reflection of
$\Omega^{\bf k}_{q,m,p}({\bf B})$ into the lower half-plane.

The domain $\Omega$ is symmetric with respect to the real axis and is convex in the direction of
imaginary axis which means that an intersection of every vertical straight line with $\Omega$ is
either empty or connected. We say that a function $g$ is convex in the direction of imaginary axis
if it maps $\mathbb D$ onto a domain convex in the direction of imaginary axis. The integral
representation of functions convex in the direction of imaginary axis is known, see, e.g.,
\cite{Rakhmanov}, \cite{Prokhorov}. Namely, such function $g$ is given by the formula
\begin{equation}
g(\zeta)=C_1\int_0^{\zeta}\frac{h(t)dt}{(1-\overline\zeta_1t)(1-\overline\zeta_2t)}+C_2,
\label{int}
\end{equation}
with constant real numbers $C_1$ and $C_2$, and a function $h$, $h(0)=1$, which is holomorphic in
$\mathbb D$ and satisfies the condition $=\re(e^{i\gamma}h(\zeta))>0$, $\zeta\in\mathbb D$,
$\gamma$ is fixed, $\gamma\in(-\pi/2,\pi/2)$, while the points $\zeta_1,\zeta_2$,
$|\zeta_1|=|\zeta_2|=1$, are disposed so that all boundary points of the image of $\mathbb D$ under
the mapping $w=\zeta[(1-\overline\zeta_1\zeta)(1-\overline\zeta_2\zeta)]^{-1}$ belong to the
straight line $\{w:\;\re(e^{i\gamma}w)=0\}$. Since the constructed function $f$ is convex in the
direction of imaginary axis and $f$ has real Taylor coefficients, we conclude that $f$ is
represented by integral (\ref{int}) with $\gamma=0$, $\zeta_1=-1$, $\zeta_2=1$. From the other
side, every function (\ref{int}) is convex in the direction of imaginary axis, and it has real
Taylor coefficients if $\gamma=0$, $\zeta_1=-1$, $\zeta_2=1$.

The aim of the proof is to show that a function $h$ in representation (\ref{int}) for $f$ is
uniquely determined by the boundary conditions.

The function $h$ in (\ref{int}) is represented by the Herglotz integral, see, e.g., \cite[p.22]{Duren},
\begin{equation}
h(\zeta)=\int_{-\pi}^{\pi}\frac{1+\zeta e^{-i\varphi}}{1-\zeta e^{-i\varphi}}d\mu(\varphi),
\;\;\;|\zeta|<1, \label{Her}
\end{equation}
where $d\mu$ is a positive measure, $\mu(\pi)-\mu(-\pi)=1$. The geometry of the domain $\Omega$
determines a special behavior of $\mu$. Denote $I:=(\alpha_{p+1},\alpha_0)$ and
$-I:=(-\alpha_0,-\alpha_{p+1})$. The function $\mu$ has to be piecewise constant on
$[-\pi,\pi]\setminus((-I)\cup I)$ with jumps at $\varphi=-\alpha_j$ and $\varphi=\alpha_j$,
$j=0,1,\dots,p+1$. Therefore,
\begin{equation}
h(\zeta)=\sum_{j=0}^{p+1}\lambda_j\left(\frac{1+\zeta e^{i\alpha_j}}{1-\zeta e^{i\alpha_j}}
+\frac{1+\zeta e^{-i\alpha_j}}{1-\zeta e^{-i\alpha_j}}\right)+\int_{(-I)\cup I} \frac{1+\zeta
e^{-i\varphi}}{1-\zeta e^{-i\varphi}}d\mu(\varphi), \;\;\;|\zeta|<1. \label{jum}
\end{equation}
Here $\mu(\varphi)$ is continuous on $(-I)\cup I$, $d\mu(\varphi)=d\mu(-\varphi)$, $\lambda_j>0$,
$j=0,1,\dots,p+1$, and
\begin{equation}
2\sum_{j=0}^{p+1}\lambda_j+2\int_Id\mu(\varphi)=1. \label{uni}
\end{equation}

Positive numbers $\lambda_0,\lambda_1,\dots,\lambda_{p+1},C_1$ and $C_2$ and the measure $d\mu$ on
$I$ and $(-I)$ in (\ref{Her}) and (\ref{jum}) are the accessory parameters of the problem. We will
show that a function (\ref{int}) with a suitable choice of the accessory parameters will coincide
with the needed function $f(\zeta)=\phi(z(\zeta))$.

Choose $C_2$ to satisfy the condition $f(-1)=0$ which, according to (\ref{int}) gives the
representation
\begin{equation}
f(\zeta)=C_1\int_{-1}^{\zeta}\frac{h(t)dt}{(1-\overline\zeta_1t)(1-\overline\zeta_2t)}.
\label{int2}
\end{equation}

The function $f$ in representations (\ref{int2}) and (\ref{jum}) has singular points at
$e^{i\alpha_0}$, $e^{i\alpha_1},\dots,$ $e^{i\alpha_p}$, $e^{i\alpha_{p+1}}$, and $e^{-i\alpha_0}$,
$e^{-i\alpha_1},\dots,e^{-i\alpha_p}$, $e^{-i\alpha_{p+1}}$. Note that representation (\ref{int2})
does not imply singularities at $\zeta=-1$ and $\zeta=1$ because $h(-1)=h(1)=0$.

Examine the character of singularity at $\zeta=e^{i\alpha_q}$. In a neighborhood of
$\zeta=e^{i\alpha_q}$ the function (\ref{int2}) has an expansion
\begin{equation}
f(\zeta)=\frac{C_1\lambda_q}{i\sin\alpha_q}\log(1-\zeta e^{-i\alpha_q})+H_q(\zeta), \label{juq}
\end{equation}
where $H_q$ is holomorphic in a neighborhood of $\zeta=e^{i\alpha_q}$. The argument of $(1-\zeta
e^{-i\alpha_q})$ has a jump equal to $\pi$ when $\zeta$ moves along the unit circle through the
point $e^{i\alpha_q}$.

Thus, as $\zeta$ moves along the unit circle through $e^{i\alpha_q}$, $\re f(\zeta)$ has a jump
$d_q$,
\begin{equation}
d_q=\frac{C_1\lambda_q\pi}{\sin\alpha_q}. \label{diq}
\end{equation}
Evidently, $d_q$ is equal to a distance between the imaginary axis and a vertical ray which is the
image of the arc $\{\zeta=e^{i\theta}:\alpha_{q-1}<\theta<\alpha_q\}$ under $w=f(\zeta)$.

Similarly, in a neighborhood of $\zeta=e^{i\alpha_{q-1}}$, $f(\zeta)$ has an expansion
\begin{equation}
f(\zeta)=\frac{C_1\lambda_{q-1}}{i\sin\alpha_{q-1}}\log(1-\zeta e^{-i\alpha_{q-1}})+
H_{q-1}(\zeta), \label{jq1}
\end{equation}
where $H_{q-1}$ is holomorphic in a neighborhood of $\zeta=e^{i\alpha_{q-1}}$. As $\zeta$ moves
along the unit circle through $e^{i\alpha_{q-1}}$, $\re f(\zeta)$ has a jump $d_{q-1}$,
\begin{equation}
d_{q-1}=\frac{C_1\lambda_{q-1}\pi}{\sin\alpha_{q-1}}, \label{dq1}
\end{equation}
where $d_{q-1}$ is equal to a distance between the vertical rays which are the images of the arcs
$\{\zeta=e^{i\theta}:\alpha_{q-1}<\theta<\alpha_q\}$ and
$\{\zeta=e^{i\theta}:\alpha_{q-2}<\theta<\alpha_{q-1}\}$ under $w=f(\zeta)$.

As for the domain $\Omega$, the distance between the imaginary axis and $l_q$ in $\Omega$ equals
$\pi k_q$, and the distance between $l_q$ and $l_{q-1}$ in $\Omega$ equals $\pi k_{q-1}$. In order
to preserve a proportion between the distances in $\Omega$ and the corresponding distances in
$f(\mathbb D)$, according to (\ref{diq}) and (\ref{dq1}), we have
\begin{equation}
\lambda_q=\lambda_{q-1}\frac{k_q}{k_{q-1}}\frac{\sin\alpha_q}{\sin\alpha_{q-1}}. \label{pr1}
\end{equation}

In the same way,
$$\lambda_{q-1}=\lambda_{q-2}\frac{k_{q-1}}{k_{q-2}}\;\frac{\sin\alpha_{q-1}}{\sin\alpha_{q-2}},
\;\dots,\;\lambda_2=\lambda_1\frac{k_2}{k_1}\;\frac{\sin\alpha_2}{\sin\alpha_1},$$ which implies
that
\begin{equation}
\lambda_j=\lambda_1\frac{k_j}{k_1}\;\frac{\sin\alpha_j}{\sin\alpha_1},\;\dots,\;j=2,\dots,q.
\label{la1}
\end{equation}

Besides, since a distance between $l_1$ and the curvilinear boundary of $\Omega$ equals
$\pi(k_0-{1\over2})$, we have
\begin{equation}
\lambda_0=\lambda_1\frac{k_0-{1\over2}}{k_1}\;\frac{\sin\alpha_0}{\sin\alpha_1}. \label{la0}
\end{equation}

From the other side,
$$\lambda_{q+1}=\lambda_{q+2}\frac{k_{q+1}}{k_{q+2}}\;\frac{\sin\alpha_{q+1}}{\sin\alpha_{q+2}},\;
\dots,\;\lambda_{p-1}=\lambda_p\frac{k_{p-1}}{k_p}\;\frac{\sin\alpha_{p-1}}{\sin\alpha_p},$$ which
implies that
\begin{equation}
\lambda_j=\lambda_p\frac{k_j}{k_p}\;\frac{\sin\alpha_j}{\sin\alpha_p},\;j=q+1,\dots,p-1.
\label{laj}
\end{equation}

Besides, since a distance between $l_p$ and the curvilinear boundary of $\Omega$ equals
$\pi(m-{1\over2})$, we have
\begin{equation}
\lambda_{p+1}=\lambda_p\frac{m-{1\over2}}{k_p}\;\frac{\sin\alpha_{p+1}}{\sin\alpha_p}. \label{lam}
\end{equation}

To compare $\lambda_1$ and $\lambda_p$, write one more relation
\begin{equation}
\lambda_p=\lambda_1\frac{k_p}{k_1}\;\frac{\sin\alpha_p}{\sin\alpha_1}. \label{lap}
\end{equation}

Take into account (\ref{la1}) - (\ref{lap}) and observe that there remain three undetermined
accessory parameters $C_1$, $\lambda_1$ and $d\mu$ on $I$ and $(-I)$ subject to relation
(\ref{uni}). Choose the scaling parameter $C_1$ so that the distance between the images of
$\{\zeta=e^{i\theta}:\alpha_1<\theta<\alpha_2\}$ and
$\{\zeta=e^{i\theta}:\alpha_0<\theta<\alpha_1\}$ under $w=f(\zeta)$ has to be equal to the distance
between $l_1$ and $l_2$,
\begin{equation}
C_1=\frac{k_1\sin\alpha_1}{\lambda_1}. \label{vC1}
\end{equation}

The positive measure $d\mu$ on $(-I)\cup I$ is of more complicated nature than the remaining number
accessory parameter $\lambda_1\in(0,1/2)$. The variation $\mu(\alpha_0-0)-\mu(\alpha_{p+1}+0)$
depends on $\lambda_1$ by (\ref{uni}),
$$\text{var}_I\mu:=\int_Id\mu(\varphi)=\frac{1}{2}-\sum_{j=0}^{p+1}\lambda_j,$$ where
$\lambda_0,\lambda_2,\dots,\lambda_{p+1}$ are linear functions of $\lambda_1$.

Define $d\mu$ on $(-I)\cup I$. The integral in the right-hand side of (\ref{jum}) can be considered
as the Schwarz integral with a smooth density which is equal to the real part of the mapping from
$(-I)\cup I$. Require that $\mu(\varphi)$ has a continuous derivative $\mu'(\varphi)$ on $(-I)\cup I$
that represents the density of the Schwarz integral. So $$\int_{(-I)\cup I}\frac{1+\zeta
e^{-i\varphi}}{1-\zeta e^{-i\varphi}}d\mu(\varphi)=\int_{(-I)\cup I}\mu'(\varphi)\frac{1+\zeta
e^{-i\varphi}}{1-\zeta e^{-i\varphi}}d\varphi,\;\;\;|\zeta|<1,$$ and
\begin{equation}\label{eq4}
2\pi\mu'(\varphi)=\re h(e^{i\varphi}),\;\;\;\varphi\in I.
\end{equation}

Let $\varphi\in I$ and $$f(e^{i\varphi})=u(\varphi)+iv(\varphi)$$ be given by (\ref{int2}). The
value $u(\alpha_{p+1})-u(\alpha_0)$ determines the length of a projection of $f(I)$ on the real
axis. Require that this length is equal to $\pi$, i.e.,
$$u(\alpha_{p+1})-u(\alpha_0)=C_1\re\int_{e^{i\alpha_0}}^{e^{i\alpha_{p+1}}}\frac{h(t)dt}{1-t^2}=
C_1\re\int_{\alpha_0}^{\alpha_{p+1}}\frac{h(e^{i\varphi})ie^{i\varphi}d\varphi}{1-e^{i2\varphi}}=$$
$$\frac{C_1}{2}\int_{\alpha_{p+1}}^{\alpha_0}\re\frac{h(e^{i\varphi})d\varphi}{\sin\varphi}=
C_1\pi\int_{\alpha_{p+1}}^{\alpha_0}\frac{\mu'(\varphi)d\varphi}{\sin\varphi}=\pi.$$ So consider
$d\mu$ satisfying the restriction
\begin{equation}
C_1\int_{\alpha_{p+1}}^{\alpha_0}\frac{d\mu(\varphi)}{\sin\varphi}=1. \label{res1}
\end{equation}

Evaluate $u(\alpha)$, $\alpha\in I$,
\begin{equation}
u(\alpha)=C_1\re\int_{-1}^{e^{i\alpha}}\frac{h(t)dt}{1-t^2}=C_1\re\left[\int_{-1}^{e^{i\alpha_0}}
\frac{h(t)dt}{1-t^2}+\int_{e^{i\alpha_0}}^{e^{i\alpha}}\frac{h(t)dt}{1-t^2}\right]= \label{reI}
\end{equation}
$$\pi\left(\sum_{j=0}^qk_j-\frac{1}{2}\right)+C_1\pi\int_{\alpha}^{\alpha_0}
\frac{\mu'(\varphi)d\varphi}{\sin\varphi}.$$

Now evaluate $v(\alpha)$, $\alpha\in I$,
\begin{equation}
v(\alpha)=C_1\im\int_{-1}^{e^{i\alpha}}\frac{h(t)dt}{1-t^2}=C_1\im\left[\int_{-1}^0
\frac{h(t)dt}{1-t^2}+\int_0^{e^{i\alpha}}\frac{h(t)dt}{1-t^2}\right]= \label{imI}
\end{equation}
$$C_1\im\int_0^{e^{i\alpha}}\frac{h(t)dt}{1-t^2}=C_1\left[\sum_{j=0}^{p+1}
\frac{\lambda_j}{\sin\alpha_j}\log\frac{\sin\frac{\alpha+\alpha_j}{2}}
{\sin\frac{|\alpha-\alpha_j|}{2}}+\right.$$ $$\left.\int_{\alpha_{p+1}}^{\alpha_0}\frac{\mu'(\varphi)}{\sin\varphi}
\log\frac{\sin\frac{\alpha+\varphi}{2}}{\sin\frac{|\alpha-\varphi|}{2}}d\varphi\right].$$

The curvilinear part of $\partial\Omega$ consists of the curve $L$ in the upper half-plane and its
reflection $\overline L$. The function $\mu'$, if it exists, provides the needed boundary behavior
of function $f$ on $(-I)\cup I$ so that $f$ maps $(-I)\cup I$ onto $L\cup\overline L$. This is
equivalent to the relation between $u(\alpha)$ and $v(\alpha)$, $\alpha\in I$, in the form of the
equation
\begin{equation}
v(\alpha)=\text{arccosh}\left(\frac{\text{cosh}B_0}{\cos(u(\alpha)-\pi\sum_{j=0}^qk_j)}\right),
\;\;\;\alpha\in I, \label{cur1}
\end{equation}
for a suitable number $B_0>0$. Equation (\ref{cur1}) is the explicit representation of $L$. Take
into account (\ref{reI}) and (\ref{imI}) to rewrite (\ref{cur1}) in the form
\begin{equation}
C_1\left[\sum_{j=0}^{p+1} \frac{\lambda_j}{\sin\alpha_j}\log\frac{\sin\frac{\alpha+\alpha_j}{2}}
{\sin\frac{|\alpha-\alpha_j|}{2}}+\int_{\alpha_{p+1}}^{\alpha_0}\frac{\mu'(\varphi)}{\sin\varphi}
\log\frac{\sin\frac{\alpha+\varphi}{2}}{\sin\frac{|\alpha-\varphi|}{2}}d\varphi\right]= \label{cur}
\end{equation}
$$\text{arccosh}\left(\frac{\text{cosh}B_0}{\sin(C_1\pi\int_{\alpha}^{\alpha_0}
\frac{\mu'(\varphi)}{\sin\varphi}d\varphi)}\right),\;\;\;\alpha\in I,$$ subject to restrictions
(\ref{uni}) and (\ref{res1}).

Show that restrictions (\ref{uni}) and (\ref{res1}) give a non-degenerate interval for values of
$\lambda_1$. Note that, due to (\ref{la1}) - (\ref{lap}),
\begin{equation}
\sum_{j=0}^{p+1}\lambda_j=\lambda_1A, \label{equl}
\end{equation}
where
\begin{equation}
A=\frac{1}{k_1\sin\alpha_1}\left[
\sum_{j=1}^pk_j\sin\alpha_j+k_0-{1\over2}+\left(m-{1\over2}\right) \sin\alpha_{p+1}\right].
\label{equA}
\end{equation}
The restrictions (\ref{uni}) and (\ref{res1}) imply that
$$\frac{k_1\sin\alpha_1}{2(1+Ak_1\sin\alpha_1)}<\lambda_1<\frac{k_1\sin\alpha_1}
{2(\sin\alpha_{p+1}+Ak_1\sin\alpha_1)},\;\;\;A>1.$$

Equation (\ref{cur}) is an integral nonlinear equation with respect to $\mu'(\varphi)$, $\varphi\in
I$. The value $B_0$ suits restriction (\ref{res1}). For a given $\mu'$, relation (\ref{uni})
determines the value of the accessory parameter $\lambda_1$.

We have to prove now that there exists an integrable solution $\mu'$ to equation (\ref{cur}). We
will approximate $\mu$ on $I$ by piecewise constant functions. To this purpose, for every natural
$n\geq1$, consider a partition $\{\beta_{1,n},\dots,\beta_{n,n}\}$ of $I$,
$$\beta_{j,n}=\alpha_{p+1}+\frac{\alpha_0-\alpha_{p+1}}{n+1}j,\;\;\;j=1,\dots,n,$$
$$\Delta_{j,n}=(\beta_{j-1,n},\beta_{j,n}),\;\;j=1,\dots,n+1,\;\;\beta_{0,n}:=\alpha_{p+1},\;\;
\beta_{n+1,n}:=\alpha_0.$$

Denote by $\mu^{(n)}$ a piecewise constant function on $I$ which has $n$ discontinuity points at
$\beta_{1,n},\dots,\beta_{n,n}$ and $\mu=\mu^{(n)}+\sum_{j=q+1}^{p+1}\lambda_j$ in a neighborhood
of $\alpha_{p+1}$. Let us construct an algorithm to define jump values $\mu_{k,n}>0$ of $\mu^{(n)}$
at the points $\beta_{k,n}$, $k=1,\dots,n$,
$$C_1\sum_{k=1}^n\frac{\mu_{k,n}}{\sin\beta_{k,n}}=1.$$

According to (\ref{jum}) with $d\mu$ substituted by $d\mu^{(n)}$, the measure $d\mu^{(n)}$
generates a function $h_n$. Let $f_n$ be defined by (\ref{int2}) with $h$ substituted by $h_n$.
Then $f_n$ maps $\mathbb D$ onto a domain $\Omega_n$ which can be described geometrically in the
following way. The domain $\Omega$ includes the part $E\cup\overline E$ of the vertical strip $S$
of width $\pi$ bounded by $L$ and $\overline L$, where $\overline E$ is a reflection of $E$.
Substitute $E\cup\overline E$ by the strip $S$ slit along $n$ vertical rays $\mathcal
L_{1,n},\dots,\mathcal L_{n,n}$ in the upper half-plane and their reflections and obtain
$\Omega_n$.

Denote by $w_{1,n},\dots,w_{n,n}$ the tips of rays $\mathcal L_{1,n},\dots,\mathcal L_{n,n}$,
respectively,
$$\re w_{k,n}=\pi\left(\sum_{j=0}^qk_j+{1\over2}\right)-\sum_{j=1}^k
\frac{C_1\pi\mu_{k,n}}{\sin\beta_{k,n}}, \;\;\;k=1,\dots,n.$$

The points $w_{k,n}$ are images of $e^{i\varphi_{k,n}}$, $k=1,\dots,n$, under the map $f_n$. All
values $\varphi_{1,n},\dots,\varphi_{n,n}$ depend on $\mu_{1,n},\dots,\mu_{n,n}$. It follows from
results of \cite{Zakharov} up to the M\"obius transformation that, for every $k=1,\dots,n$,
$\varphi_{k,n}(\mu_{1,n},\dots,\mu_{n,n})$ tends asymptotically to the center point of
$\Delta_{k,n}$ as $\mu_{k,n}\to0$.

The function $u_n(\varphi)=\re f_n(e^{i\varphi})$, $\varphi\in I$, is piecewise constant,
$$u_n(\varphi)=\re w_{k,n},\;\;\varphi\in\Delta_{k,n},\;\;k=1,\dots,n.$$

The function $v_n(\varphi)=\im f_n(e^{i\varphi})$ is continuous on all $\Delta_{k,n}$,
$k=1,\dots,n$,
$$v_n(\varphi)=C_1\left[\sum_{j=0}^{p+1}
\frac{\lambda_j}{\sin\alpha_j}\log\frac{\sin\frac{\varphi+\alpha_j}{2}}
{\sin\frac{|\varphi-\alpha_j|}{2}}+\sum_{j=1}^n\frac{\mu_{j,n}}{\sin\beta_{j,n}}
\log\frac{\sin\frac{\varphi+\beta_{j,n}}{2}}{\sin\frac{|\varphi-\beta_{j,n}|}{2}}\right]$$ for
$\varphi\in\Delta_{k,n}$, $k=1,\dots,n$, $$v_n(\varphi_{k,n})=\im
w_{k,n}=\min_{\varphi\in\Delta_{k,n}}v_n(\varphi),\;\;k=1,\dots,n.$$

We have $n$ free parameters $\mu_{1,n},\dots,\mu_{n,n}$ subject to restriction (\ref{res1}) and
consider $B_0$ as one more free parameter $B_{0,n}$ to establish $n$ values $v_n(\varphi_{k,n})$,
$k=1,\dots,n$, satisfying (\ref{cur1}). If $B_{0,n}$ tends to 0, then $\Omega_n$ degenerates. If
$B_{0,n}$ tends to infinity, then $E$ degenerates. Values $v_n(\varphi_{k,n})$ are strictly
monotone in $\mu_{j,n}$, $j,k=1,\dots,n$, and $B_{0,n}$.

For $n>1$, call $(\mu_{1,n},\dots,\mu_{n,n})$ {\it admissible} if
$$C_1\sum_{k=1}^{n-1}\frac{\mu_{k,n}}{\sin\beta_{k,n}}<1.$$ Vary $B_{0,n}$ to find $\mu_{n,n}$ for
which restriction (\ref{res1}) and condition (\ref{cur1}) are satisfied at $\alpha=\varphi_{n,n}$
with $u=u_n(\varphi_{n,n})=\re w_{n,n}$ and $v=v_n(\varphi_{n,n})=\im w_{n,n}$. In this way, we
find $\mu_{n,n}=\mu_{n,n}(\mu_{1,n},\dots,\mu_{n-1,n})$ and
$B_{0,n}=B_{0,n}(\mu_{1,n},\dots,\mu_{n-1,n})$ for any admissible $(\mu_{1,n},\dots,\mu_{n-1,n})$.

Next, for $n>2$ and admissible $(\mu_{1,n},\dots,\mu_{n-2,n})$ such that
$$C_1\sum_{k=1}^{n-2}\frac{\mu_{k,n}}{\sin\beta_{k,n}}<1,$$ find
$\mu_{n-1,n}=\mu_{n-1,n}(\mu_{1,n},\dots,\mu_{n-2,n})$ for which (\ref{cur1}) is satisfied at
$\alpha=\varphi_{n-1,n}$ with $u=u_n(\varphi_{n-1,n})=\re w_{n-1,n}$ and
$v=v_n(\varphi_{n-1,n})=\im w_{n-1,n}$ where $\mu_{n,n}=\mu_{n,n}(\mu_{1,n},\dots,\mu_{n-1,n})$ and
$B_{0,n}=B_{0,n}(\mu_{1,n},\dots,\mu_{n-1,n})$. A solution exists since $\mu_{n-1,n}$ depends
continuously on $\mu_{1,n},\dots,\mu_{n-2,n}$ and the two extreme positions of admissible
$\mu_{n-1,n}$ lead to degenerate configurations.

Continue this process up to the last step and find all $\mu_{1,n},\dots,\mu_{n,n}$ and $B_{0,n}$
solving (\ref{cur1}) at every $\alpha=\varphi_{k,n}$, $k=1,\dots,n$, with $u=u_n$ and $v=v_n$.

When $n$ tends to infinity, Helly's theorems admit existence of a convergent subsequence from
$\{\mu^{(n)}\}$ for which $B_{0,n}$ tends to a certain $B_0$. The corresponding subsequence of
functions $h_n$ generated by $\{d\mu^{(n)}\}$ converges to a function $h$ and the subsequence
$\{f_n\}$ given by (\ref{int2}) with $h_n$ in its right-hand side converges to $f^*$. Preserve the
denotation $\mu^{(n)}$ for the convergent subsequence. The variation for the limit measure
satisfies (\ref{res1}). The limit function $f^*$ is not constant because of the normalization
condition (\ref{res1}). Therefore, $f^*$ maps $\mathbb D$ onto a kernel $\Omega^*$ of the
subsequence of domains $\{\Omega_n\}$. Evidently, $\Omega^*=\Omega$ if the rays $\{\mathcal
L_{k,n}\}$, $1\leq k\leq n$, $n\geq1$, are dense in the strip $S$. This is equivalent to the
condition
$$\lim_{n\to\infty}\max_{1\leq k\leq n}\mu_{k,n}=0.$$

Suppose that there exists a sequence $\{k_n\}$, $1\leq k\leq n$, (for a suitable subsequence of $n$'s) such that $$
\lim_{n\to\infty}\mu_{k_n,n}=\mu_0>0.$$ This means that $\partial\Omega^*$ contains a vertical
straight segment $\mathcal L_0$ in the closure of the strip $S$. A preimage of $\mathcal L_0$ under
$f^*$ is an arc in $\partial\mathbb D$. Let $e^{i\varphi^*}$ be an inner point of this arc,
$\varphi^*\in I$. In a neighborhood of $\varphi^*$ there is a set of $\Delta_{k,n}$ with $n$ large
enough and certain integers $k$, $w_{k,n}\in\Delta_{k,n}$. By construction,
$$\frac{\im(w_{k,n}-w_{j,n})}{\re(w_{k,n}-w_{j,n})}$$ is bounded. To the contrary, this ratio is
unbounded if $\Delta_{k,n}$ and $\Delta_{j,n}$ are in a neighborhood of $\varphi^*$. The
contradiction implies that the supposition is wrong, and the rays $\{\mathcal L_{k,n}\}$, $1\leq
k\leq n$, $n\geq1$, are dense in $S$.

Finally, we obtained a smooth measure $d\mu$ which is the limit of subsequence of $\{d\mu^{(n)}\}$.
Equation (\ref{uni}) with this $d\mu$ gives us $\sum_{j=0}^{p+1}\lambda_j$. Hence, according to
(\ref{equl}) and (\ref{equA}), we determine the last free parameter $\lambda_1$ and complete the
proof.

\textit{Proof of Theorem 2.}
Proof is quite similar to the proof of \cite[Theorem 4]{ereyud1} (compare also \cite[Theorem 3.1]{nazetal}).
First of all we note that $$f=1-(-1)^{\sum_{j=0}^qk_j}L\cos\phi$$
is real on the positive real semi-axis and pure imaginary on the positive imaginary semi-axis. So by two reflections with respect to the coordinate axes, $f$ extends to an odd function analytic in $ \mathbb{C}\backslash\{0,\pm x_1,\ldots,\pm x_p\}. $ The region $\Omega^{\bf
k}_{q,m,p}({\bf B}^*)$ is close to the strip
$$\left\{w:\mathrm{Re} w\in\left(\pi\sum_{j=1}^qk_j,\pi\left(\sum_{j=0}^qk_j-1/2\right)\right)\right\}$$
as Im $w\to\infty$ reached along the right-hand side of $l_1,$ and to the strip $$\left\{w:\mathrm{Re} w\in\left(\pi\left(\sum_{j=0}^qk_j+1/2\right),\pi\left(\sum_{j=0}^qk_j+m\right)\right)\right\}$$
as Im $w\to\infty$ reached along the left-hand side of $l_p.$ So $\phi\sim(2k_0-1)\log1/z,$ $z\to 0,$ $\phi\sim(2m-1)\log z,$ $z\to\infty,$ $\phi\sim k_j\log(1/(z\pm x_j)),$ $z\to\pm x_j,$ $j=1,\ldots,p,$ and, therefore $f$ is a rational function of the form (\ref{form}). Next we note that the graph of $f$ alternates $\sum_{j=0}^pk_j+m+1$ times on $[a,1]$ between $1-L$ and $1+L.$ Now we conclude the proof by using the general Chebyshev alternation theorem on the uniform approximation of continuous functions by Chebyshev systems (here we consider the Chebyshev system $$\frac{1}{x^{2k_0-1}(x^2-x_1^2)^{k_1}\ldots(x^2-x_p^2)^{k_p}},\ldots,\frac{x^n}{x^{2k_0-1}(x^2-x_1^2)^{k_1}\ldots(x^2-x_p^2)^{k_p}}$$ on $[-1,a]\cup[a,1].)$ 

Finally we note there is a unique point on the imaginary axis where $\phi=\pi\sum_{j=0}^qk_j+iB_0^*,$ what gives (\ref{deviation}).


\begin{thebibliography}{4}

\bibitem{akh} Akhiezer N.I. \"Uber einige Funktionen, die in gegebenen Intervallen am wenigsten von Null abweichen,  Bull. Soc. Phys.-Mathem. Kazan. Ser.3,  3. N2, 1-69 (1928).

\bibitem{apt} Aptekarev, A.I.
Asymptotic properties of polynomials orthogonal on a system of contours, and periodic motions of Toda lattices, Math. USSR, Sb. 53, 233-260 (1986).

\bibitem{bog} Bogatyrev A., Extremal polynomials and Riemann surfaces, Springer, Berlin, Heidelberg (2012).

\bibitem{SIAM1} Diakonikolas I., Gopalan P., Jaiswal R., Servedio R.A., Viola E., Bounded independence fools halfspaces, SIAM J. Comput., 39, 3441-3462 (2010).

\bibitem{Duren}
Duren P.L.,  Univalent Functions, Springer, Heidelberg, New York (1983).

\bibitem{ereyud1} Eremenko A., Yuditskii P., Uniform approximation of sgn $x$ by polynomials and entire functions, J. Anal. Math., 101, 313-324 (2007).

\bibitem{ereyudsurv} Eremenko A., Yuditskii P., Comb functions, in: Recent advances in orthogonal polynomials, special functions, and their applications (J. Arves\'u, G. L\'opez Lagomasino, Eds.), AMS: Providence,RI, 99-118 (2012).

\bibitem{leb} Lebedev V.I., Extremal polynomials and methods of optimization of numerical algorithms, Sbornik: Math., 195, 1413-1460 (2004).

\bibitem{luk1} Lukashov A.L., On Chebyshev-Markov rational fractions over several intervals, J. Approx.
Theory, 95, 333-352 (1998).

\bibitem{luk2} Lukashov A. L. Inequalities for the derivatives of rational functions on several intervals, Izv. Math. 68, 543-565 (2004).

\bibitem{lukpeh} Lukashov A. L., Peherstorfer F. Automorphic orthogonal and extremal polynomials, Canad. J. Math., 55, 576-608 (2003).

\bibitem{luktys}   Lukashov A.L.,  Tyshkevich S.V., Extremal rational functions on several arcs of the unit circle with zeros on these arcs, Izv. Saratov. Univ. Mat. Mekh. Inform., 9,no.1 , 8-13 (2009) [in Russian].
  
\bibitem{nazetal} Nazarov F., Peherstorfer F., Volberg A., Yuditskii P., Asymptotics of the best polynomial approximation of $|x|^p$ and of the best Laurent polynomial approximation of sgn$(x)$ on two symmetric intervals, Constr. Approx., 29, 23-39 (2009).

\bibitem{peh} Peherstorfer F. Orthogonal and extremal 
polynomials on several intervals,  J. Comp. Appl. Math.,48, 187-205 (1993). 
 
\bibitem{ereyud2} Peherstorfer F., Yuditskii P. Uniform approximation of sgn $x$ by rational functions with prescribed poles, J. Mathem. Phys., Analysis, Geom., 3, 95-108 (2007).

\bibitem{Prokhorov}
Prokhorov D.V.,  Level curves of functions convex in the direction of an axis, Math. Notes,  44, 767-769 (1988).

\bibitem{Rakhmanov}
Prokhorov D.V., Rakhmanov B.N., Integral representations of a class of univalent functions,
Math. Notes 19, 24-28 (1976).

\bibitem{Zakharov}
Prokhorov D., Zakharov A., Harmonic measures of sides of a slit perpendicular to the domain
boundary, J. Math. Anal. Appl. 394, 738-743 (2012).

\bibitem{rob}  Robinson R.M. Conjugate algebraic integers in real point sets, Math. Zeit.,84, 415-427 (1964).

\bibitem{SIAM2} Sherstov A.A., Strong direct product theorems for quantum communication and query complexity, SIAM J. Comput., 41, 1122-1165 (2012).

\bibitem{sodyud} Sodin M.L., Yuditskij P.M., Functions deviating least from zero on closed subsets of the real axis, St. Petersbg. Math. J., 4, 201-249 (1993).

\bibitem{wid} Widom H.,  Extremal polynomials associated with a system of curves in the complex plane, Adv. Math., 3, 127-232 (1969).

\end{thebibliography}
\end{document}